\documentclass[a4paper]{amsart}

\RequirePackage{amsmath}
\RequirePackage{bm}
\RequirePackage{amssymb}
\RequirePackage{upref}
\RequirePackage{amsthm}
\RequirePackage{enumerate}
\RequirePackage{pb-diagram}
\RequirePackage{amsfonts}
\RequirePackage[mathscr]{eucal}
\RequirePackage{verbatim}
\RequirePackage{xr}
\RequirePackage{graphicx}
\usepackage{calc}
\usepackage{xspace}
\RequirePackage{color}
\RequirePackage{ifthen}

\newcommand{\al}{\alpha}
\newcommand{\bet}{\beta}
\newcommand{\ga}{\gamma}
\newcommand{\de}{\delta }
\newcommand{\e}{\epsilon}

\newcommand{\f}{\varphi}
\newcommand{\h}{\eta}

\newcommand{\tht}{\theta}
\newcommand{\ka}{\kappa}
\newcommand{\lam}{\lambda}
\newcommand{\m}{\mu}
\newcommand{\n}{\nu}

\newcommand{\vt}{\vartheta}

\newcommand{\s}{\sigma}
\newcommand{\x}{\xi}

\newcommand{\C}{\varGamma}

\newcommand{\Om}{\varOmega}

\newcommand{\Si}{\varSigma}

\newcommand{\nbdd}{\nobreakdash--}

\newcommand{\fmo}[1]{F_{|_{#1}}}
\newcommand{\fu}[3]{#1\hspace{0pt}_{|_{#2_#3}}}
\newcommand{\fv}[2]{#1\hspace{0pt}_{|_{#2}}}

\newcommand{\so}{{\mc S_0}}

\newcommand{\const}{\tup{const}}

\newcommand{\msp[1]}[1]{\mspace{#1mu}}
\newcommand{\low}[1]{{\hbox{}_{#1}}}

\newcommand{\R}[1][n+1]{{\protect\mathbb R}^{#1}}

\newcommand{\N}{{\protect\mathbb N}}

\newcommand{\eR}{\stackrel{\lower1ex \hbox{\rule{6.5pt}{0.5pt}}}{\msp[3]\R[]}}
\newcommand{\eN}{\stackrel{\lower1ex \hbox{\rule{6.5pt}{0.5pt}}}{\msp[1]\N}}
\newcommand{\eO}{\stackrel{\lower1ex \hbox{\rule{6pt}{0.5pt}}}{\msc O}}

\DeclareMathOperator{\graph}{graph}

\newcommand\pa{\partial}
\newcommand\pde[2]{\frac {\partial#1}{\partial#2}}
\newcommand\pd[3]{\frac {\partial#1}{\partial#2^#3}}   
\newcommand\pdc[3]{\frac {\partial#1}{\partial#2_#3}}   
\newcommand\pdm[4]{\frac {\partial#1}{\partial#2_#3^#4}}   
\newcommand\pddc[4]{\frac {{\partial\hskip0.15em}^2#1}{\partial {#2_
#3}\,\partial{#2_#4}}} 

\newcommand{\un}{\infty}
\newcommand{\A}{\forall}

\newcommand{\set}[2]{\{\,#1\colon #2\,\}}
\newcommand{\uu}{\cup}
\newcommand{\ii}{\cap}
\newcommand{\uuu}{\bigcup}

\newcommand{\uud}{ \stackrel{\lower 1ex \hbox {.}}{\uu}}
\newcommand{\uuud}[1]{ \stackrel{\lower 1ex \hbox {.}}{\uuu_{#1}}}
\newcommand\su{\subset}

\newcommand{\sminus}[1][28]{\raise 0.#1ex\hbox{$\scriptstyle\setminus$}}

\newcommand{\abs}[1]{\lvert#1\rvert}

\newcommand{\norm}[1]{\lVert#1\rVert}

\newcommand{\nnorm}[1]{| \mspace{-2mu} |\mspace{-2mu}|#1| \mspace{-2mu}
|\mspace{-2mu}|}
\newcommand{\spd}[2]{\protect\langle #1,#2\protect\rangle}

\newcommand\cha[3]{{\bar\varGamma}_{#1#2}^#3}

\newcommand{\riema}[4]{{\bar R}_{#1#2#3#4}}

\newcommand{\tbf}{\textbf}
\newcommand{\tit}{\textit}

\newcommand{\tup}{\textup}

\newcommand{\mc}{\protect\mathcal}
\newcommand{\msc}{\protect\mathscr}

\providecommand{\bysame}{\makebox[3em]{\hrulefill}\thinspace}

\newcommand{\ci}{\cite}
\newcommand{\bib}{\bibitem}

\newcommand{\bt}{\begin{thm}}
\newcommand{\bl}{\begin{lem}}
\newcommand{\bc}{\begin{cor}}
\newcommand{\bd}{\begin{definition}}
\newcommand{\bpp}{\begin{prop}}
\newcommand{\br}{\begin{rem}}
\newcommand{\bn}{\begin{note}}
\newcommand{\be}{\begin{ex}}
\newcommand{\bes}{\begin{exs}}
\newcommand{\bb}{\begin{example}}
\newcommand{\bbs}{\begin{examples}}
\newcommand{\ba}{\begin{axiom}}
\newcommand{\bas}{\begin{assumption}}

\newcommand{\et}{\end{thm}}
\newcommand{\el}{\end{lem}}
\newcommand{\ec}{\end{cor}}
\newcommand{\ed}{\end{definition}}
\newcommand{\epp}{\end{prop}}
\newcommand{\er}{\end{rem}}
\newcommand{\en}{\end{note}}
\newcommand{\ee}{\end{ex}}
\newcommand{\ees}{\end{exs}}
\newcommand{\eb}{\end{example}}
\newcommand{\ebs}{\end{examples}}
\newcommand{\ea}{\end{axiom}}
\newcommand{\eas}{\end{assumption}}

\newcommand{\bp}{\begin{proof}}
\newcommand{\ep}{\end{proof}}
\newcommand{\eps}{\renewcommand{\qed}{}\end{proof}}

\newcommand{\bal}{\begin{align}}

\newcommand{\bi}[1][1.]{\begin{enumerate}[\upshape #1]}
\newcommand{\bia}[1][(1)]{\begin{enumerate}[\upshape #1]}
\newcommand{\bin}[1][1]{\begin{enumerate}[\upshape\bfseries #1]}
\newcommand{\bir}[1][(i)]{\begin{enumerate}[\upshape #1]}
\newcommand{\bic}[1][(i)]{\begin{enumerate}[\upshape\hspace{2\cma}#1]}
\newcommand{\bis}[2][1.]{\begin{enumerate}[\upshape\hspace{#2\parindent}#1]}
\newcommand{\ei}{\end{enumerate}}

\newcommand\ndots{\raise 0.47ex \hbox {,}\hskip0.06em\cdots %
     \raise 0.47ex \hbox {,}\hskip0.06em} 

\newcommand{\q}{\quad}
\newcommand{\qq}{\qquad}

\newcommand{\hp}{\hphantom}

\newcommand\nd{\noindent}

\newskip\Csmallskipamount                                                
\Csmallskipamount=\smallskipamount
\newskip\Cmedskipamount
\Cmedskipamount=\medskipamount
\newskip\Cbigskipamount
\Cbigskipamount=\bigskipamount

\newcommand\cvs{\vspace\Csmallskipamount}   
\newcommand\cvm{\vspace\Cmedskipamount}
\newcommand\cvb{\vspace\Cbigskipamount}

\newskip\csa
\csa=\smallskipamount

\newskip\cma
\cma=\medskipamount

\newskip\cba
\cba=\bigskipamount

\newdimen\spt
\newcommand\citem{\cvs\advance\itemno by
1{(\romannumeral\the\itemno})\hskip3pt}
\newcommand{\bitem}{\cvm\nd\advance\itemno by
1{\bf\the\itemno}\hspace{\cma}}

\newcount\itemno
\newcommand{\las}[1]{\label{S:#1}}

\newcommand{\lae}[1]{\label{E:#1}}
\newcommand{\lat}[1]{\label{T:#1}}
\newcommand{\lal}[1]{\label{L:#1}}

\newcommand{\lap}[1]{\label{P:#1}}
\newcommand{\lar}[1]{\label{R:#1}}

\newcommand{\rs}[1]{Section~\ref{S:#1}}

\newcommand{\rl}[1]{Lemma~\ref{L:#1}}

\newcommand{\rp}[1]{Proposition~\ref{P:#1}}
\newcommand{\rr}[1]{Remark~\ref{R:#1}}

\newcommand{\re}[1]{\eqref{E:#1}}

\newskip\thmskip
\thmskip=\parindent

\newskip\hsk
\newenvironment{hinw}{\labelsep=0pt\begin{list}{}{\labelsep=0pt\itemindent=0pt\labelwidth=0pt\leftmargin=\parindent\rightmargin=0pt\partopsep=\cba}%
\item\it\nopagebreak\nopagebreak}%
{\end{list}}

\newcommand\bh{\begin{hinw}}
\newcommand{\eh}{\end{hinw}}

\newtheoremstyle{normal}
  {\cba}
  {\cba}
  {}
  {\thmskip}
  {\bfseries}
  {.}
  {\hsk}
  {}

\newtheoremstyle{abschnitt}
  {\cba}
  {\cba}
  {}
  {\thmskip}
  {\bfseries}
  {.}
  {\hsk}
  {}

\newtheoremstyle{italic}
  {\cba}
  {\cba}
  {\itshape}
  {\thmskip}
  {\bfseries}
  {.}
  {\hsk}
  {}

\newtheoremstyle{aufgaben}
  {\cba}
  {\cba}
  {}
  {}
  {\normalsize\bfseries}
  {.}
  {\hsk}
  {}

\newtheoremstyle{break}
  {\cba}
  {\cba}
  {\itshape}
  {}
  {\bfseries}
  {.}
  {\newline}
  {}

\swapnumbers
\theoremstyle{italic}
\newtheorem{thm}[subsection]{Theorem}
\newtheorem{lem}[subsection]{Lemma}
\newtheorem{prop}[subsection]{Proposition}
\newtheorem{cor}[subsection]{Corollary}

\theoremstyle{normal}
\newtheorem{rem}[subsection]{Remark}
\newtheorem{definition}[subsection]{Definition}
\newtheorem{example}[subsection]{Example}
\newtheorem{examples}[subsection]{Examples}
\newtheorem{ex}[subsection]{Exercise}
\newtheorem{note}[subsection]{}
\newtheorem{axiom}[subsection]{Axiom}
\newtheorem{assumption}[subsection]{Assumption}

\theoremstyle{aufgaben}
\newtheorem{exs}[subsection]{Exercises}

\swapnumbers

\numberwithin{equation}{section}
\numberwithin{figure}{section}

\newenvironment{textequation}[1][0.8]
{\begin{equation}
\begin{aligned}
\begin{minipage}{#1\linewidth}}
{\end{minipage}
\end{aligned}
\end{equation}
\ignorespacesafterend}

\newcommand{\btext}{\begin{textequation}}
\newcommand{\etext}{\end{textequation}}

\def\hinweis{\@startsection{subsection}{2}%
 \z@{0.7\linespacing\@plus 0.5\linespacing}{0.7\linespacing}%
{\normalfont\itshape\indent}}

\newcounter{hours}\newcounter{minutes}
\newcommand{\printtime}{%
\setcounter{hours}{\time/60}%
\setcounter{minutes}{\time-\value{hours}*60}%
\ifthenelse{\value{minutes}<9}{\thehours :0\theminutes}{\thehours:\theminutes}}

\usepackage[german,english]{babel}
\usepackage{graphicx}
\RequirePackage{amsmath}
\RequirePackage{bm}
\RequirePackage{amssymb}
\RequirePackage{upref}
\RequirePackage{amsthm}
\RequirePackage{enumerate}
\RequirePackage{pb-diagram}
\RequirePackage{amsfonts}
\RequirePackage[mathscr]{eucal}
\RequirePackage{verbatim}
\RequirePackage{xr}
\RequirePackage{graphicx}
\usepackage{calc}
\usepackage{xspace}

\makeatletter
\RequirePackage{color}
\newcommand{\ann}[1]{\renewcommand{\@makefnmark}{\mbox{$^{\color{red}{\@thefnmark}}$}}%
\footnote {#1}}
\makeatother

\RequirePackage{upref}
\RequirePackage{amsthm}
\RequirePackage{enumerate}
\usepackage[mathscr]{eucal}

\usepackage{xr-hyper}

\usepackage{calc}

\newlength{\oddsidemarginlength}
\newlength{\topmarginlength}

\newcounter{numberoflines}
\newcounter{tempcc}
\oddsidemargin=\oddsidemarginlength
\evensidemargin=\oddsidemargin
\topmargin=\topmarginlength
\usepackage[colorlinks=true,linkcolor=blue,citecolor=blue,urlcolor=blue]{hyperref}

\begin{document}

\title[]{The scalar curvature flow in Lorentzian manifolds}

\author{Christian Enz}
\email{c.enz@stud.uni-heidelberg.de}

\begin{abstract}
We prove the existence of closed hypersurfaces of prescribed scalar curvature in globally hyperbolic Lorentzian manifolds provided there are barriers.
\end{abstract}

\maketitle

\tableofcontents

\setcounter{section}{-1}
\section{Introduction}\las{0}

\cvb\nd
We want to find a closed spacelike hypersurface of prescribed scalar curvature in a globally hyperbolic Lorentzian manifold $N$ having a compact Cauchy hypersurface $\mc S_0$. Looking at the Gau{\ss} equation for a spacelike hypersurface $M$, we deduce that its scalar curvature $R$ satisfies
\begin{equation}\lae{0.1}
R=-[H^2-|A|^2]+\bar{R}+2\bar{R}_{\alpha\beta}\nu^\alpha \nu^\beta.
\end{equation}
Denoting the curvature operator defined by $H_2$ by $F$, then this equation is equivalent to
\begin{equation}\lae{0.2}
R=-2F(h_{ij})+\bar R+2\bar R_{\al\bet}\n^{\al}\n^{\bet}.
\end{equation}
Thus, we have to allow that the right-hand side $f$ of the equation
\begin{equation}\lae{0.3}
F_{|_M}=f
\end{equation}
\nd is defined in $T(N)$, or more precisely, after choosing a local trivialization of $T(N)$, that $f$ depends on $x\in N$ and timelike vectors $\n\in T_x(N)$, and we look for a closed spacelike hypersurface $M$ satisfying
\begin{equation}\lae{0.4}
F_{|_M}=f(x,\n)\qq\A x\in M.
\end{equation}

In \ci{cg10} Gerhardt solved this problem by using the method of elliptic regularization. We give a new existence proof, based on the curvature estimates in \ci{cg13}, by showing that the scalar curvature flow exists for all time, and that the leaves $M(t)$ of the flow converge to a solution of \re{0.4}.

To give a precise statement of the existence result we need a few definitions and
assumptions. First, we assume that $\Om$ is a precompact, connected, 
open subset of $N$, that is bounded by  two compact, connected, spacelike
hypersurfaces $M_1$ and $M_2$ of class $C^{6,\al}$, where $M_1$ is supposed
to lie in the past of
$M_2$.

Let $F=H_2$ be the scalar curvature operator defined on the open cone
$\C_2\su \R[n]$, and
$f=f(x,\n)$ be of class
$C^{4,\al}$ in its arguments such that 
\begin{align}
0<c_1&\le f(x,\n)\qq\tup{if}\q\spd\n\n=-1,\lae{0.5}\\[\cma]
\nnorm{f_\bet(x,\n)}&\le c_2 (1+\nnorm\n^2),\lae{0.6}\\
\intertext{and}
\nnorm{f_{\n^\bet}(x,\n)}&\le c_3 (1+\nnorm\n),\lae{0.7}
\end{align}
for all $x\in\bar\Om$ and all past directed timelike vectors $\n\in T_x(\Om)$,
where $\nnorm{\cdot}$ is a Riemannian reference metric that will be detailed in \rs{1}.

\cvb
\br
The condition \re{0.5} is reasonable as is evident from the Einstein equation
\begin{equation}\lae{0.8}
\bar R_{\al\bet}-\tfrac12 \bar R \msp \bar g_{\al\bet}=T_{\al\bet},
\end{equation} 
where the energy-momentum tensor $T_{\al\bet}$ is supposed to be positive
semi-definite for timelike vectors (weak energy condition, cf. \ci[p.~89]{HE});
but it would be
convenient, if the estimate in \re{0.5}
would be valid for all timelike vectors. In fact, we may assume this without loss
of generality: Let $\vt$ be a smooth real function such that
\begin{equation}\lae{0.9}
\frac{c_1}2\le \vt\q\tup{and}\q \vt(t)=t\q \A\,t\ge c_1,
\end{equation}
\nd then, we can replace $f$ by $\vt\circ f$ and the new function satisfies our
requirements for all timelike vectors. We therefore assume in the following that the relation \re{0.5} holds for all
timelike vectors $\n\in T_x(N)$ and all $x\in \bar\Om$.
\er

\cvb
We suppose that the boundary components $M_i$ act as barriers for  $(F,f)$.

\nd $M_2$ is an upper barrier for $(F,f)$, if $M_2$ is admissible, i.e. its
principal curvatures $(\ka_i)$ with respect to the past directed normal belong to
$\C_2$, and if
\begin{equation}\lae{0.10}
\fv F{M_2}\ge f(x,\n)\qq\A\,x\in M_2.
\end{equation}
$M_1$ is a lower barrier for $(F,f)$, if at the points $\Si\su M_1$, where $M_1$
is admissible, there holds
\begin{equation}\lae{0.11}
\fv F\Si \le f(x,\n)\qq\A\,x\in \Si.
\end{equation}
\nd $\Si$ may be empty.

Now, we can state the main theorem.

\cvb
\bt\lat{0.2}
Let $M_1$ be a lower and $M_2$ an upper barrier for $(F,f)$, where $F=H_2$. Then, the problem
\begin{equation}
\fmo M= f(x,\n)
\end{equation}
has an admissible solution $M\su \bar\Om$ of class $C^{6,\al}$ that can be
written as a graph over $\mc S_0$ provided there exists a strictly convex
function $\chi\in C^2(\bar\Om)$.
\et

\cvb
\br
As proved in \ci[Lemma 2.7]{cg8} the existence of a strictly convex
function
$\chi$ is guaranteed by the assumption that the level hypersurfaces
$\{x^0=\tup{const}\}$ are strictly convex in $\bar\Om$, where $(x^\al)$ is a
Gaussian coordinate system associated with $\so$.

Looking at  Robertson-Walker space-times it seems that the assumption of the
existence of  a strictly convex function in the neighbourhood of a given compact
set is not too restrictive: in Minkowski space e.g. $\chi=-\abs{x^0}^2 +\abs x^2$
is a globally defined strictly convex function.
\er
\cvb\nd

\section{Notations and preliminary results}\las{1}

\cvb \nd We refer to \ci{cg10} or to \ci{cg12} for a more detailed treatment.

\cvb

\section{Curvature functions}\las{2}

\cvb\nd
Let $\C\su\R[n]$ be an open  cone containing the positive cone $\C_+$, and
$F\in C^{2,\al}(\C)\ii C^0(\bar\C)$ a positive symmetric function satisfying the
condition
\begin{equation}\lae{2.1}
F_i=\pd F\ka i>0\; ,
\end{equation}
then, $F$ can also be viewed as a function defined on the space of symmetric
 matrices $\msc S_\Gamma$, the eigenvalues of which belong to $\C $, namely, let
$(h_{ij})\in
\msc S_\Gamma$ with eigenvalues $\ka_i,\,1\le i\le n$, then define $F$ on $\msc S_\Gamma$
by
\begin{equation}\lae{2.2}
F(h_{ij})=F(\ka_i).
\end{equation}
If we define
\begin{align}
F^{ij}&=\pde F{h_{ij}}\\
\intertext{and}
F^{ij,kl}&=\pddc Fh{{ij}}{{kl}}
\end{align}
then,
\begin{equation}\lae{2.5}
F^{ij}\x_i\x_j=\pdc F\ka i \abs{\x^i}^2\q\A\, \x\in\R[n],
\end{equation}
in an appropriate coordinate system,
\begin{equation}\lae{2.6}
F^{ij} \,\text{is diagonal if $h_{ij}$ is diagonal,}
\end{equation}
and
\begin{equation}\lae{2.7}
F^{ij,kl}\h_{ij}\h_{kl}=\pddc F{\ka}ij\h_{ii}\h_{jj}+\sum_{i\ne
j}\frac{F_i-F_j}{\ka_i-\ka_j}(\h_{ij})^2,
\end{equation}
for any $(\h_{ij})\in \msc S$, where $\msc S$ is the space of all symmetric matrices. The second term on the right-hand side of \re{2.7} is non-positive if
$F$ is concave, and non-negative if $F$ is convex, and has to be interpreted as a
limit if $\ka_i=\ka_j$.

The preceding considerations are
also applicable if the
$\ka_i$ are the principal curvatures of a spacelike hypersurface $M$ with metric
$(g_{ij})$.
$F$ can then be looked at as being defined on the space of all symmetric tensors
$(h_{ij})$ the eigenvalues of which 
belong to
$\C$. Such tensors will be called admissible; when the second
fundamental form of $M$ is admissible, then, we also call $M$ admissible.

For an admissible tensor $(h_{ij})$
\begin{equation}\lae{2.8}
F^{ij}=\pdc Fh{{ij}}
\end{equation}
is  a contravariant tensor of second order. Sometimes it will be convenient
to circumvent the dependence on the metric by considering $F$ to depend on the
mixed tensor
\begin{equation}\lae{2.9}
h_j^i=g^{ik}h_{kj}.
\end{equation}
Then,
\begin{equation}\lae{2.10}
F_i^j=\pdm Fhji
\end{equation}
is also a mixed tensor with contravariant index $j$ and covariant index $i$. Such functions $F$ are called curvature functions.

Important examples are the elementary
symmetric polynomials of order $k$, $H_k$, $1\le k\le n$,
\begin{equation}\lae{2.11}
H_k(\ka_i)=\sum_{i_1<\cdots <i_k} \ka_{i_1}\cdots \ka_{i_k}.
\end{equation}
They are defined on an open set $\C_k$  that can be
characterized as the connected component of $\{H_k>0\}$ that contains the positive cone
$\C_+$.
The $\C_k$ are cones, $\C_n=\C_+$, and in \ci{lg} it is proved that
\begin{equation}\lae{2.12}
\C_k\su\C_{k-1}.
\end{equation}
Huisken and Sinestrari in \ci[Section 2]{hs} gave an equivalent characterisation of $\C_k$ by showing that
\begin{equation}\lae{2.13}
\C_k=\{(\ka_i)\in \R[n]: H_1(\ka_i)>0, H_2(\ka_i)>0, \dots , H_k(\ka_i)>0\}.
\end{equation}
They also proved that $\C_k$ is convex.
The $H_k$ are strictly monotone in $\C_k$, cf. \ci[Lemma 2.4]{hs}, and the the k-th roots
\begin{equation}\lae{2.14}
\s_k=H_k^\frac1k
\end{equation}
are also concave, cf. \ci{dm}.

Since we have in mind that the $\ka_i$ are the principal curvatures of a
hypersurface, we use the standard symbols $H$ and $\abs A$ for
\begin{align}
H&=\sum_i\ka_i,\\
\intertext{and}
\abs A^2&=\sum_i \ka^2_i.
\end{align}
We note that
\begin{equation}\lae{2.17}
\tfrac1nH^2\le\abs A^2.
\end{equation}
The scalar curvature function $F=H_2$ can be expressed as
\begin{equation}\lae{2.18}
F=\tfrac12{(H^2-\abs A^2)},
\end{equation}
and we deduce that for $(\ka_i)\in \C_2$
\begin{align}
\abs A^2&\le H^2,\lae{2.19}\\[\cma]
F&\le \tfrac12 H^2,\lae{2.20}\\[\cma]
F_i&=H-\ka_i,\lae{2.21}\\
\intertext{and hence,}
H&>\ka_i,\lae{2.22}\\[\cma]
H F_i&\ge F,\lae{2.23}
\end{align}
for \re{2.23} is equivalent to
\begin{equation}\lae{2.24}
H\ka_i\le \tfrac12 H^2 + \tfrac12 \abs A^2,
\end{equation}
which is obviously valid.

\cvb
\section{The evolution problem}\las{3}

\cvb\nd
To prove the existence of hypersurfaces of prescribed curvature $F$ for $F=\s_2$ we look at the evolution problem
\begin{equation}
\begin{aligned}\lae{3.1}
\dot x&=(F-f)\n,\\[\cma]
x(0)&=x_0,
\end{aligned}
\end{equation}
where $\n$ is the past-directed normal of the flow hypersurfaces $M(t)$,
$F=\s_2$ the curvature evaluated at $M(t)$,  $x=x(t)$ an embedding and $x_0$ an
embedding of an initial hypersurface $M_0$, which we choose to be the upper
barrier $M_2$.

Since $F$ is an elliptic operator,  short-time existence, and hence, existence in a
maximal time interval $[0,T^*)$ is guaranteed, cf. \ci{cg12}. If we are able to prove uniform a priori estimates in $C^{2,\al}$, long-time existence and convergence to a
stationary solution will follow immediately.

But before we prove the a priori estimates, we want to show how the metric, the
second fundamental form, and the normal vector of the hypersurfaces $M(t)$
evolve. All time derivatives are total derivatives. We shall omit the proofs, which can be found in \ci[Chapter 2.3]{cg12}.

\cvb
\bl
The metric, the normal vector and the second fundamental form of $M(t)$ satisfy the evolution equations
\begin{align}
\dot g_{ij}&=2( F- f)h_{ij},\lae{3.2}\\[\cma]
\dot \n&=\nabla_M( F- f)=g^{ij}( F- f)_i x_j,\lae{3.3}\\[\cma]
\dot h_i^j&=( F- f)_i^j- ( F- f) h_i^k h_k^j-( F- f) \riema
\al\bet\ga\de\n^\al x_i^\bet \n^\ga x_k^\de g^{kj},\lae{3.4}\\[\cma]
\dot h_{ij}&=(F-f)_{ij}+( F- f)h_i^k h_{kj}-(F-f)\riema
\al\bet\ga\de\n^\al x_i^\bet \n^\ga x_j^\de.\lae{3.5}
\end{align}
\el
\cvb
\bl\lal{3.2}
The term $( F- f)$ evolves according to the equation
\begin{multline}\lae{3.6}
{( F- f)}^\prime- F^{ij}( F- f)_{ij}=\msp[3]- 
F^{ij}h_{ik}h_j^k ( F- f)
 - f_\al\n^\al ( F- f)\\-f_{\n^\al}x^\al_i(F-f)_jg^{ij}
- F^{ij}\riema \al\bet\ga\de\n^\al x_i^\bet \n^\ga x_j^\de ( F- f).
\end{multline}
\el
\cvb\nd
From \re{3.4} we deduce with the help of the Ricci identities a parabolic equation
for the second fundamental form, cf. \ci[Lemma 2.4.9]{cg12}.
\cvb
\bl\lal{3.3}
The mixed tensor $h_i^j$ satisfies the parabolic equation
\begin{equation}\lae{3.7}
\begin{aligned}
\dot h_i^j- &F^{kl}h_{i;kl}^j\\
&=-  F^{kl}h_{rk}h_l^rh_i^j
+f h_i^kh_k^j\\
&\hp{+}- f_{\al\bet} x_i^\al x_k^\bet g^{kj}-  f_\al\n^\al
h_i^j-f_{\al\n^\bet}(x^\al_i x^\bet_kh^{kj}+x^\al_l
x^\bet_k h^{k}_i\, g^{lj})\\ &\hp{=}
-f_{\n^\al\n^\bet}x^\al_lx^\bet_kh^k_ih^{lj}-f_{\n^\bet} x^\bet_k h^k_{i;l}\,g^{lj} 
-f_{\n^\al}\n^\al h^k_i h^j_k\\ &\hp{=}+
F^{kl,rs}h_{kl;i}h_{rs;}^{\hphantom{rs;}j}+2 F^{kl}\riema
\al\bet\ga\de x_m^\al x_i ^\bet x_k^\ga x_r^\de h_l^m g^{rj}\\
&\hp{=}- F^{kl}\riema \al\bet\ga\de x_m^\al x_k ^\bet x_r^\ga x_l^\de
h_i^m g^{rj}- F^{kl}\riema \al\bet\ga\de x_m^\al x_k ^\bet x_i^\ga x_l^\de h^{mj} \\
&\hp{=}- F^{kl}\riema \al\bet\ga\de\n^\al x_k^\bet\n^\ga x_l^\de h_i^j+ f
\riema \al\bet\ga\de\n^\al x_i^\bet\n^\ga x_m^\de g^{mj}\\
&\hp{=}+ F^{kl}\bar R_{\al\bet\ga\de;\e}\{\n^\al x_k^\bet x_l^\ga x_i^\de
x_m^\e g^{mj}+\n^\al x_i^\bet x_k^\ga x_m^\de x_l^\e g^{mj}\}.
\end{aligned}
\end{equation}
\el
\cvb
\br\lar{3.4}
In view of the maximum principle, we immediately deduce from \re{3.6} that the
term $( F- f)$ has a sign during the evolution if it has one at the beginning,
i.e., if the starting hypersurface $M_0$ is the upper barrier $M_2$, then
$( F- f)$ is non-negative
\begin{equation}\lae{3.8}
F\ge f.
\end{equation}
\er
\cvb
\section{Lower order estimates}\las{4}

\cvb\nd
Since the two boundary components $M_1, M_2$ of $\pa\Om$ are compact, connected
spacelike hypersurfaces, they can be written as graphs over the Cauchy
hypersurface $\so$, $M_i=\graph u_i$, $i=1,2$, and we have
\begin{equation}\lae{4.1}
u_1\le u_2,
\end{equation}
for $M_1$ should lie in the past of $M_2$.

Let us look at the evolution equation \re{3.1} with initial hypersurface $M_0$
equal to $M_2$ defined on  a maximal time interval
$I=[0,T^*)$, $T^*\le
\un$. 
Since the initial hypersurface is a graph over $\so$, we can write
\begin{equation}\lae{4.2}
M(t)=\graph\fu{u(t)}S0\q \A\,t\in I,
\end{equation}
where $u$ is defined in the cylinder $Q_{T^*}=I\times \so$.

We then deduce from
\re{3.1}, looking at the component $\al=0$, that $u$ satisfies a parabolic
equation of the form
\begin{equation}\lae{4.3}
\dot u=-e^{-\psi}v^{-1}(F- f),
\end{equation}
where we  use the notations in \rs{2}, and where we emphasize that the time
derivative is a total derivative, i.e.
\begin{equation}\lae{4.4}
\dot u=\pde ut+u_i\dot x^i.
\end{equation}
Since the past directed normal can be expressed as
\begin{equation}\lae{4.5}
(\n^\al)=-e^{-\psi}v^{-1}(1,u^i),
\end{equation}
we conclude from \re{3.1}, \re{4.3} and \re{4.4}
\begin{equation}\lae{4.6}
\pde ut=-e^{-\psi}v(F- f).
\end{equation}
Thus, $\pde ut$ is non-positive in view of \rr{3.4}.

Next, let us state our first a priori estimate, \ci[Theorem 2.7.9]{cg12}.
\cvb
\bl\lal{4.1}
Suppose that the boundary components act as barriers for $(F,f)$, then  the flow hypersurfaces stay in
$\bar
\Om$ during
the evolution.
\el
\cvb
\nd For the $C^1$-estimate the term $\tilde v=v^{-1}$ is of great importance. It satisfies the following evolution equation.
\cvb
\bl\lal{4.2}
Consider the flow \re{3.1} in the distinguished coordinate system associated
with $\so$. Then, $\tilde v$ satisfies the evolution equation
\begin{equation}\lae{4.7}
\begin{aligned}
\dot{\tilde v}- F^{ij}\tilde v_{ij}=&- F^{ij}h_{ik}h_j^k\tilde v
-f\h_{\al\bet}\n^\al\n^\bet\\
&-2 F^{ij}h_j^k x_i^\al x_k^\bet \h_{\al\bet}- F^{ij}\h_{\al\bet\ga}x_i^\bet
x_j^\ga\n^\al\\
&- F^{ij}\riema \al\bet\ga\de\n^\al x_i^\bet x_k^\ga x_j^\de\h_\e x_l^\e g^{kl}\\
&- f_\bet x_i^\bet x_k^\al \h_\al g^{ik} -f_{\n^\bet}x^\bet_k h^{ik}x^\al_i\h_\al,
\end{aligned}
\end{equation}
where $\h$ is the covariant vector field $(\h_\al)=e^{\psi}(-1,0,\dotsc,0)$.
\el

\cvb\nd
The proof uses the relation
\begin{equation}\lae{4.8}
\tilde v=\h_\al \n^\al
\end{equation}
and
is identical to that of \ci[Lemma 4.4]{cg8} having in
mind that  presently $f$  also depends on $\n$.
\cvb
\bl\lal{4.3}
Let $M(t)=\graph u(t)$ be the flow hypersurfaces, then, we have
\begin{equation}\lae{4.9}
\begin{aligned}
\dot u-F^{ij}u_{ij}=e^{-\psi}\tilde v f+\cha 000\, F^{ij}u_i u_j 
+2\msp
F^{ij}\cha 0i0\,u_j +F^{ij}\cha ij0,
\end{aligned}
\end{equation}
where all covariant derivatives are taken with respect to the induced metric of
the flow hypersurfaces, and the time derivative $\dot u$ is the total time
derivative, i.e. it is given by \re{4.4}.
\el
\cvb
\bp
We use the relation \re{4.3}.
\ep
\cvb\nd
As an immediate consequence we obtain
\cvb
\bl\lal{4.4}
The composite function
\begin{equation}\lae{4.10}
\f=e^{\m e^{\lam u}},
\end{equation}
where $\m,\lam$ are constants, satisfies the equation
\begin{equation}\lae{4.11}
\begin{aligned}
\dot\f -F^{ij}\f_{ij}=&fe^{-\psi} \tilde v\msp \m \lam \msp[2] e^{\lam u} \,\f + F^{ij}
u_i u_j \,\cha000\,\m \lam\,e^{\lam u}\f\\[\cma]
&+2\msp F^{ij} u_i \cha 0j0\,\m \lam\msp[2] e^{\lam u}\,\f+ F^{ij}\cha
ij0\,\m\lam\msp[2] e^{\lam u}\,\f\\[\cma]
&-[1+\m\msp[2] e^{\lam u}] F^{ij} u_i u_j\,\m \lam^2\msp[2] e^{\lam u}\,\f.
\end{aligned}
\end{equation}
\el
\cvb\nd
Before we can prove the $C^1$- estimates we need two more lemmata.
\cvb
\bl\lal{4.5}
There is a constant $c=c(\Om)$ such that for any positive function $0<\e=\e(x)$
on $\so$ and any hypersurface $M(t)$ of the flow we have
\begin{align}
\nnorm{\n}&\le c \msp \tilde v,\lae{4.12}\\[\cma]
g^{ij}&\le c\msp \tilde v^2\msp \s^{ij},\lae{4.13}\\[\cma]
F^{ij}&\le F^{kl} g_{kl}\msp g^{ij},\lae{4.14}
\end{align}

\begin{equation}\lae{4.15}
\begin{aligned}
\abs{F^{ij}h^k_j x^\al_i x^\bet_k \msp\h_{\al\bet}}\le \frac\e 2F^{ij}h^k_i
h_{kj}\msp
\tilde v + \frac c{2\e} F^{ij}g_{ij}\msp  \tilde v^3,
\end{aligned}
\end{equation}
\begin{equation}\lae{4.16}
\abs{F^{ij}\h_{\al\bet\ga} x^\bet_i x^\ga_j \n^\al}\le c\msp \tilde v^3 F^{ij}g_{ij},
\end{equation}
\begin{equation}\lae{4.17}
\abs{F^{ij}\riema \al\bet\ga\de \n^\al x^\bet_i x^\ga_k x^\de_j\h_\e x^\e_l g^{kl}}\le
c\msp
\tilde v^3 F^{ij}g_{ij}.
\end{equation}
\el

\cvb
\bl\lal{4.6}
Let $M\su \bar \Om$ be a graph over $\so$, $M=\graph u$, and
$\e=\e(x)$ a function defined on $\so$, $0<\e <\frac12$. Let $\f$ be
defined through
\begin{equation}\lae{4.18}
\f=e^{\m e^{\lam u}},
\end{equation}
where $0<\m$ and $\lam<0$. Then, there exists $c=c(\Om)$ such that
\begin{equation}\lae{4.19}
\begin{aligned}
2\abs{F^{ij} \tilde v_i \f_j}&\le c\msp F^{ij}g_{ij} \tilde v^3 \abs{\lam} \m e^{\lam u}
\f +(1-2\e) F^{ij} h^k_i h_{kj} \tilde v \f\\[\cma]
&\hp{\le} +\frac1{1-2\e} F^{ij} u_i u_j \m^2 \lam^2 e^{2\lam u} \tilde v \f.
\end{aligned}
\end{equation}
\el

\cvb\nd
A proof of \rl{4.5} and \rl{4.6} can be found in \ci{cg10}.

Applying \rl{4.5} to the evolution equation for $ \tilde v$ we conclude
\cvb
\bl\lal{4.7}
There exists a constant $c=c(\Om)$ such that for any function $\e$,
$0<\e=\e(x)<1$, defined on $\so$ the term $ \tilde v$ satisfies an evolution
inequality of the form
\begin{equation}\lae{4.20}
\begin{aligned}
\dot {\tilde v} -F^{ij} \tilde v_{ij}&\le -(1-\e) F^{ij} h^k_i h_{kj} \tilde v -f
\h_{\al\bet} \n^\al \n^\bet\\[\cma]
&\hp{\le} +\frac{c}\e F^{ij} g_{ij} \tilde v^3 +c\msp \nnorm{f_\bet} \tilde v^2+
f_{\n^\bet} x^\bet_l h^{kl} u_k e^\psi.
\end{aligned}
\end{equation}
\el
\cvb\nd
We are now ready to prove the uniform boundedness of $ \tilde v$.
\cvb
\bpp\lap{4.8}
Assume that there are positive constants $c_i$, $1\le i\le 3$, such that for any
$x\in \Om$ and any past directed timelike vector $\n$ there holds
\begin{align}
-c_1&\le f(x,\n),\lae{4.21}\\[\cma]
\nnorm{f_\bet(x,\n)}&\le c_2 (1+\nnorm{\n}),\lae{4.22}\\
\intertext{and}
\nnorm{f_{\n^\bet}(x,\n)}&\le c_3.\lae{4.23}
\end{align}
Then, the term $ \tilde v$ remains uniformly bounded during the evolution
\begin{equation}\lae{4.24}
\tilde v\le c=c(\Om, c_1, c_2, c_3).
\end{equation}
\epp
\cvb
\bp
We show that the function
\begin{equation}\lae{4.25}
w= \tilde v\f,
\end{equation}
$\f$ as in \re{4.18}, is uniformly bounded, if we choose
\begin{equation}\lae{4.26}
0<\m<1\q \tup{and}\q \lam<<-1,
\end{equation}
appropriately,  and assume furthermore, without loss of generality, that $u\le
-1$, for otherwise replace $u$ by $(u-c)$, $c$ large, in the definition of $\f$.
With the help of \rl{4.4}, \rl{4.6} and \rl{4.7} we derive
from the relation
\begin{equation}\lae{4.27}
\dot w - F^{ij} w_{ij}=[\dot{ \tilde v}- F^{ij} \tilde v_{ij}] \f+
[\dot\f-F^{ij}\f_{ij}] \tilde v-2 F^{ij} \tilde v_i \f_j
\end{equation}
the parabolic inequality
\begin{equation}\lae{4.28}
\begin{aligned}
\dot w -F^{ij} w_{ij}&\le -\e\msp[2] F^{ij} h^k_i h_{kj} \tilde v \f +
c[\e^{-1}+\abs\lam \m e^{\lam u}]F^{ij} g_{ij} \tilde v^3 \f\\[\cma]
&\hp{\le} +[\frac1{1-2\e}-1] F^{ij} u_i u_j \m^2\lam^2 e^{2\lam u} \tilde v
\f \\[\cma]
&\hp{\le}-F^{ij} u_i u_j \m\lam^2 e^{\lam u} \tilde v \f \\[\cma]
&\hp{\le} +f
[-\h_{\al\bet}\n^\al \n^\bet+e^{-\psi}
\m
\lam e^{\lam u} \tilde v^2] \f\\[\cma]
&\hp{\le} + c\msp[2] \nnorm{f_\bet} \tilde v^2 \f +f_{\n^\bet} x^\bet_l h^{kl} u_k
e^\psi
\f ,
\end{aligned}
\end{equation}
where we have chosen the same function $\e=\e(x)$ in \rl{4.6} { }resp. \rl{4.7}.
We claim that $w$ is uniformly bounded  provided $\m$ and $\lam$ are chosen
appropriately.
We shall use the maximum
principle, therefore let $0<T<T^*$ and $x_0=x(t_0,\x_0)$ be such that
\begin{equation}\lae{4.29}
\sup_{[0,T]}\sup_{M(t)}w=w(t_0,\x_0).
\end{equation}
To exploit the good term
\begin{equation}\lae{4.30}
 -\e \msp[2] F^{ij} h^k_i h_{kj} \tilde v \f,
\end{equation}
we use the fact that $Dw(x_0)=0$, or, equivalently
\begin{equation}\lae{4.31}
\begin{aligned}
- \tilde v_i&= \m \lam e^{\lam u} \tilde v u_i\\[\cma]
&= e^\psi h^k_iu_k -\h_{\al\bet}\n^\al x^\bet_i,
\end{aligned}
\end{equation}
where the second equation follows from  \re{4.8} and the definition of
the covariant vectorfield
$\h=e^\psi (-1,0,\dotsc,0)$.
Next, we choose a coordinate system $(\x^i)$ such that in the critical point
\begin{equation}\lae{4.32}
g_{ij}=\de_{ij}\qq\tup{and}\qq h^k_i=\ka_i \de^k_i,
\end{equation}
and the labelling of the principal curvatures corresponds to
\begin{equation}\lae{4.33}
\ka_1\le \ka_2\le \dotsb \le \ka_n.
\end{equation}
Then, we deduce from \re{4.31}
\begin{equation}\lae{4.34}
e^\psi \ka_i u_i=\m\lam e^{\lam u} \tilde v u_i + \h_{\al\bet}\n^\al x^\bet_i.
\end{equation}
Assume that  $ \tilde v(x_0)\ge 2$, and let $i=i_0$ be an index such that
\begin{equation}\lae{4.35}
\abs{u_{i_0}}^2\ge \frac1n\norm{Du}^2.
\end{equation}
Setting $(e^i)=\pde{}{\x^{i_0}}$ and assuming without loss of generality that
$0< u_ie^i$ in $x_0$ we infer from \rl{2.6}
\begin{equation}\lae{4.36}
\begin{aligned}
e^\psi \ka_{i_0} u_ie^i&=\m\lam e^{\lam u} \tilde v u_i e^i+\h_{\al\bet}\n^\al
x^\bet_ie^i \\[\cma]
&\le \m\lam e^{\lam u} \tilde v u_ie^i +c \tilde v^2,
\end{aligned}
\end{equation}
and we deduce further in view of \re{4.35} that 
\begin{equation}\lae{4.37}
\begin{aligned}
\ka_{i_0}\le [\m\lam e^{\lam u}+ c] \tilde v e^{-\psi} 
\le \frac12\m\lam e^{\lam u} \tilde  v e^{-\psi},
\end{aligned}
\end{equation}
if $\abs\lam$ is sufficiently large,
i.e. $\ka_{i_0}$ is negative and of the same order as $ \tilde v$.
The Weingarten equation yield
\begin{equation}\lae{4.38}
\nnorm{\n^\bet_iu^i}=\nnorm{h^k_iu^i x^\bet_k}\le c \tilde v [h^k_iu^i
h_{kl}u^l]^{\frac12},
\end{equation}
and therefore, we infer from \re{4.31}
\begin{equation}\lae{4.39}
\nnorm{\n^\bet_i u^i}\le c \m\abs\lam e^{\lam u} \tilde v^3
\end{equation}
in critical points of $w$, and hence, that in those points, the term involving
$f_{\tilde \n^\bet}$ on the right-hand side of inequality \re{4.28} can be estimated
from above by
\begin{equation}\lae{4.40}
\abs{f_{  \n^\bet}\n^\bet_i  u^i e^\psi 
\f}\le c c_3 \m\abs\lam e^{\lam u} \tilde v^3\f.
\end{equation}
Next, let us estimate the crucial term in \re{4.30}. Using the particular
coordinate system \re{4.32}, as well as the inequalities \re{4.33}, together with
the fact that $\ka_{i_0}$ is negative, we conclude
\begin{equation}\lae{4.41}
-F^{ij}h^k_i h_{kj}\le -\sum_{i=1}^{i_0} F^i_i
\ka^2_i\le -\sum_{i=1}^{i_0} F^i_i \ka^2_{i_0}.
\end{equation}
$ F$ is concave, and therefore, we have in
view of \re{4.33}
\begin{equation}\lae{4.42}
F^1_1\ge F^2_2\ge\dotsb\ge  F^n_n,
\end{equation}
cf. \ci[Lemma 2]{eh}.
\cvm
Hence, we conclude
\begin{equation}\lae{4.43}
-\sum_{i=1}^{i_0}F^i_i\le -F^1_1\le -\frac1n\sum_{i=1}^nF^i_i.\\
\end{equation}
Using \re{4.37}, \re{4.41}, \re{4.43}, \re{2.20} and \re{2.21}, we
deduce further
\begin{equation}\lae{4.44}
\begin{aligned}
-F^{ij} h^k_i h_{kj}&\le-cF^{ij}g_{ij} \m^2\lam^2 e^{2\lam u} \tilde v^2\\[\cma]
&\le -c \m^2\lam^2 e^{2\lam u} \tilde v^2
\end{aligned}
\end{equation}
Inserting this estimate, and the estimate in \re{4.40} in \re{4.28}, with $\e=e^{-\lam u}$, we obtain
\begin{equation}\lae{4.45}
\begin{aligned}
\dot w-F^{ij} w_{ij}&\le-c F^{ij}g_{ij}\m^2\lam^2 e^{\lam u} \tilde v^3\f+cF^{ij}g_{ij}\m\abs\lam e^{\lam u} \tilde v^3\f
\\[\cma]
&\hp{\le}
+\frac2{1-2\e}\msp[2]F^{ij} u_i u_j \m^2\lam^2 e^{\lam u} \tilde v
\f -F^{ij} u_i u_j \m\lam^2 e^{\lam u} \tilde v \f \\[\cma]
&\hp{\le}+ c \msp c_1 \m \abs\lam e^{\lam u} \tilde v^2 \f + c\msp c_2 \tilde v^3 \f
+c\msp c_3  \m \abs\lam e^{\lam u} \tilde v^3 \f, 
\end{aligned}
\end{equation}
where $\abs\lam$ is chosen so large that
\begin{equation}\lae{4.46}
e^{-\lam u}\le \frac14.
\end{equation}
Choosing $\m=\tfrac14$ and $\abs\lam$ sufficient large, we see that in view of \re{4.44} the right-hand side of the preceding inequality is negative, contradicting the maximum principle, i.e. the
maximum of $w$ cannot occur at a point where $ \tilde v\ge 2$. Thus, the desired uniform
estimate for $w$ and hence $ \tilde v$ is proved.
\ep
\cvb
\br
Notice that the preceding $C^1$-estimate is valid for any curvature function $F$ that is monotone, concave and homogeneous of degree 1.
\er
\cvb\nd
Let us close this section with an interesting observation that is an immediate
consequence of the preceding proof, we have especially  \re{4.41} and \re{4.43} in mind.

\cvb
\bl\lal{4.9}
Suppose $F=\s_2$ is evaluated at a point $(\ka_i)$, and assume that $\ka_{i_0}$ is a
component  that is either negative or the smallest component of that particular
$n$- tupel, then
\begin{equation}\lae{4.47}
\sum_{i=1}^n F_i\ka^2_i\ge \frac1n \sum_{i=1}^n F_i \msp[2] \ka^2_{i_0}.
\end{equation}
\el

\cvb
\section{Curvature estimates}\las{5}

\cvb\nd
We want to prove that the principal curvatures of the flow hypersurfaces are
uniformly bounded.
Let us first prove an a priori estimate for $F$.
\cvb
\bl\lal{5.1}
Let $M(t)$, $0\le t<T^*$, be solutions of the evolution problem \re{3.1} with
$M(0)=M_2$ and $F=\s_2$.
Then, $F$ is bounded from above during the evolution provided the $M(t)$ are uniformly spacelike, i.e. uniform $C^1$-
estimates are valid.
\el
\cvb
\bp
Let $0<T<T^*$ and $x_0=x(t_0,\x_0)$ be such that
\begin{equation}\lae{5.1}
\sup_{[0,T]}\sup_{M(t)}(F-f)=(F-f)(x_0)> 0.
\end{equation}
Applying the maximum principle we deduce from \re{3.6}
\begin{equation}\lae{5.2}
0\le -F^{ij}h_{ik}h^{kj}+c(1+F^{ij}g_{ij}),
\end{equation}
where we have estimated bounded terms by a constant $c$.

\nd Then, we infer from \re{2.17}, \re{2.21} and \re{2.23}
\begin{equation}\lae{5.3}
0\le -\tfrac{1}{2n}FH+c(1+F^{-1}H),
\end{equation}
which is equivalent to
\begin{equation}\lae{5.4}
0\le -\tfrac{1}{2n}F^2+c(FH^{-1}+1).
\end{equation}
Thus, in view of \re{2.20}, we obtain an a priori estimate for $F$.
\ep
\cvb
\br
Let $\chi$ be the strictly convex function. Its evolution equation is
\begin{equation}\lae{5.5}
\begin{aligned}
\dot\chi- F^{ij}\chi_{ij}&=f\chi_\al\nu^\al -F^{ij}\chi_{\al\bet}x^\al_ix^\bet_j\\
&\le f\chi_\al\nu^\al -c_0  F^{ij}g_{ij},
\end{aligned}
\end{equation}
where $c_0>0$ is independent of $t$.
\er
\cvb
\bpp
Under the assumptions of Lemma 5.1 the principal curvatures $\ka_i, 1\le i\le n$, of the flow hypersurfaces are uniformly bounded during the evolution provided there exists a strictly convex function $\chi\in C^2(\bar\Om)$.
\epp
\cvb
\bp
Let $\zeta$ and $w$ be respectively defined by
\begin{align}
\zeta&=\sup\set{{h_{ij}\h^i\h^j}}{{\norm\h=1}},\lae{5.6}\\[\cma]
w&=\log\zeta+\lam \chi,\lae{5.7}
\end{align}
where $\lam>0$ is supposed to be large.  We claim that
$w$ is bounded, if $\lam$ is chosen sufficiently large. 

\nd Let $0<T<T^*$, and $x_0=x_0(t_0)$, with $ 0<t_0\le T$, be a point in $M(t_0)$ such
that
\begin{equation}\lae{5.8}
\sup_{M_0}w<\sup\set {\sup_{M(t)} w}{0<t\le T}=w(x_0).
\end{equation}
We then introduce a Riemannian normal coordinate system $(\x^i)$ at $x_0\in
M(t_0)$ such that at $x_0=x(t_0,\x_0)$ we have
\begin{equation}\lae{5.9}
g_{ij}=\delta_{ij}\q \tup{and}\q \zeta=h_n^n.
\end{equation}
Let $\tilde \h=(\tilde \h^i)$ be the contravariant vector field defined by
\begin{equation}\lae{5.10}
\tilde \h=(0,\dotsc,0,1),
\end{equation}
and set
\begin{equation}\lae{5.11}
\tilde \zeta=\frac{h_{ij}\tilde \h^i\tilde \h^j}{g_{ij}\tilde \h^i\tilde \h^j}\raise 2pt
\hbox{.}
\end{equation}
$\tilde\zeta$ is well defined in neighbourhood of $(t_0,\x_0)$.

\nd Now, define $\tilde w$ by replacing $\zeta$ by $\tilde \zeta$ in \re{5.7}; then, $\tilde w$
assumes its maximum at $(t_0,\x_0)$. Moreover, at $(t_0,\x_0)$ we have 
\begin{equation}\lae{5.12}
\dot{\tilde \zeta}=\dot h_n^n,
\end{equation}
and the spatial derivatives do also coincide; in short, at $(t_0,\x_0)$ $\tilde \zeta$
satisfies the same differential equation \re{3.7} as $h_n^n$. For the sake of
greater clarity, let us therefore treat $h_n^n$ like a scalar and pretend that $w$
is defined by 
\begin{equation}\lae{5.13}
w=\log h_n^n+\lam \chi.
\end{equation} 
We assume that the section curvatures are labelled according to \re{4.33}.

\nd At  $(t_0,\xi_0)$ we have $\dot w \ge0$, and, in view of the maximum principle, we deduce from \re{3.7}, \re{5.5}, \re{2.7} and \re{4.42}
\begin{equation}\lae{5.14}
\begin{aligned}
0&\le-\tfrac12F^{ij}h_{ki}h^k_j+c h^n_n+cF^{ij}g_{ij}+\lam c-\lam c_0F^{ij}g_{ij}\\
&\hp{\le}\;+F^{ij}(\log h^n_n)_i(\log h^n_n)_j+\frac2{\ka_n-\ka_1} \sum_{i=1}^n(F_n-F_i)(h_{ni;}^{\hp{ni;}n})^2 (h^n_n)^{-1},\\
\end{aligned}
\end{equation}
where we have estimated bounded terms by a constant $c$, and assumed that $h^n_n$ and $\lam$ are larger than $1$. We distinguish two cases

\cvb
\tit{Case} $1$.\q Suppose that
\begin{equation}\lae{5.15}
\abs{\ka_1}\ge \e_1 \ka_n,
\end{equation}
where $\e_1>0$ is small. Then, we infer from \rl{4.9}
\begin{equation}\lae{5.16}
F^{ij}h_{ki}h^k_j\ge \tfrac1n F^{ij}g_{ij}\e_1^2\ka_n^2,
\end{equation}
and 
\begin{equation}\lae{5.17}
F^{ij}g_{ij}\ge F(1,\dots,1),
\end{equation}
for a proof see \ci[Lemma 2.2.19]{cg12}.

\nd Since $Dw=0$,
\begin{equation}\lae{5.18}
D\log h^n_n=-\lam D\chi,
\end{equation}
hence
\begin{equation}\lae{5.19}
F^{ij}(\log h^n_n)_i(\log h^n_n)_j\le \lam^2F^{ij}\chi_i\chi_j.
\end{equation}
Hence, we conclude that $\ka_n$ is a priori bounded in this case.

\cvb
\tit{Case} $2$.\q Suppose that
\begin{equation}\lae{5.20}
\ka_1\ge -\e_1\ka_n,
\end{equation}
then, the last term in inequality \re{5.14} is estimated from above by
\begin{equation}\lae{5.21}
\begin{aligned}
&\frac2{1+\e_1} \sum_{i=1}^n(F_n-F_i)(h_{ni;}^{\hp{ni;}n})^2 (h^n_n)^{-2}&\\
&\le\frac2{1+2\e_1} \sum_{i=1}^n(F_n-F_i)(h_{nn;}^{\hp{nn;}i})^2 (h^n_n)^{-2}\\
&\qq+c(\e_1)\sum_{i=1}^{n}(F_i-F_n)\ka_n^{-2},
\end{aligned}
\end{equation} 
where we used the Codazzi equation. The last sum can be easily balanced.

\nd The terms in \re{5.14} containing the derivative of $h^n_n$ can therefore be estimated from above by
\begin{equation}\lae{5.22}
\begin{aligned}
&-\frac{1-2\e_1}{1+2\e_1} \sum_{i=1}^nF_i(h_{nn;}^{\hp{nn;}i})^2 (h^n_n)^{-2}\\
&+\frac2{1+2\e_1}F_n\sum_{i=1}^n(h_{nn;}^{\hp{nn;}i})^2 (h^n_n)^{-2}\\
&\le 2 F_n\sum_{i=1}^n(h_{nn;}^{\hp{nn;}i})^2 (h^n_n)^{-2}\\
&= 2\lam^2 F_n \norm{D\chi}^2.\\
\end{aligned}
\end{equation}
Hence, we infer
\begin{equation}\lae{5.23}
\begin{aligned}
0\le - \tfrac12 F_n\ka_n^2&+\lam^2c F_n+c\ka_n+cF^{ij}g_{ij} \\
&+\lam c-\lam c_0F^{ij}g_{ij}.
\end{aligned}
\end{equation}
From \re{2.21}, \re{2.22} and \rl{5.1} we deduce
\begin{equation}\lae{5.24}
F^{ij}g_{ij}\ge c\ka_n,
\end{equation}
thus, taking \re{5.17} into account, we obtain an a priori estimate 
\begin{equation}\lae{5.25}
\ka_n\le \const,
\end{equation}
if $\lam$ is chosen large enough. Notice that $\e_1$ is only subject to the requirement $0<\e_1<\frac12$.
\ep
\cvb
\br\lar{5.4}
In view of \re{0.5} and \re{3.8}, we conclude that the principal curvatures of the flow hypersurfaces stay in a compact subset of $\C$.
\er
\cvb
\section{Existence of a solution}\las{6}
\cvb\nd
We shall show that the solution of the evolution problem \re{3.1} exists for all
time, and that it converges to a stationary solution. 
\cvb
\bpp\lap{6.1}
The solutions $M(t)=\graph u(t)$ of the evolution problem \re{3.1} with $F=\s_2$ and $M(0)=M_2$ exist for all time and converge to a stationary solution provided $f\in C^{4,\al}$ satisfies the conditions \re{0.5}, \re{4.22} and
\re{4.23}.
\epp
\cvb
\bp
Let us look at the scalar
version of the flow as in
\re{4.6}
\begin{equation}\lae{6.1}
\pde ut=-e^{-\psi}v(F- f).
\end{equation}
This is  a scalar parabolic differential equation defined on the cylinder
\begin{equation}\lae{6.2}
Q_{T^*}=[0,T^*)\times \so
\end{equation}
with initial value $u(0)=u_2\in C^{4,\al}(\so)$.

In view of the a priori estimates,
which we have established in the preceding sections, we know that
\begin{equation}\lae{6.3}
{\abs u}_\low{2,0,\so}\le c
\end{equation}
and
\begin{equation}\lae{6.4}
F\,\tup{is uniformly elliptic in}\,u
\end{equation}
independently of $t$, in view of \rr{5.4}. Thus, we can apply
the known regularity results, see e.g. \ci[Chapter 5.5]{nk}, where even more
general operators are considered,  to conclude that uniform
$C^{2,\al}$-estimates are valid.
\nd
Therefore, the maximal time interval is unbounded, i.e. $T^*=\un$.

Now, integrating \re{6.1} with respect to $t$, and observing that the right-hand
side is non-positive, yields
\begin{equation}\lae{6.5}
u(0,x)-u(t,x)=\int_0^te^{-\psi}v(F- f)\ge c\int_0^t(F- f),
\end{equation}
i.e.,
\begin{equation}\lae{6.6}
\int_0^\un \abs{F- f}<\un\qq\A\msp x\in \so.
\end{equation}
Hence, for any $x\in\so$ there is a sequence $t_k\rightarrow \un$ such that
$(F- f)\rightarrow 0$.

\nd On the other hand, $u(\cdot,x)$ is monotone decreasing and therefore
\begin{equation}\lae{6.7}
\lim_{t\rightarrow \un}u(t,x)=\tilde u(x)
\end{equation}
exists and is of class $C^{2,\al}(\so)$.We
conclude that $\tilde u$ is a stationary solution, and that
\begin{equation}\lae{6.8}
\lim_{t\rightarrow \un}(F- f)=0.
\end{equation}
Now, we can deduce that uniform $C^{6,\al}$-estimates are valid, cf \ci[Theorem 6.5]{cg14}.
\nd Hence, we conclude that the functions $u(t,\cdot)$ converge in $C^6(\so)$ to $\tilde u \in C^{6,\al}(\so)$.
\ep
\cvb\nd
We want to solve the equation
\begin{equation}\lae{6.9}
\s_{2|_M}=f^{\frac12}(x,\nu),
\end{equation}
where $f$ satisfies the conditions of \re{0.5}, \re{0.6} and \re{0.7}. Thus we would like to apply the preceding existence result. But, unfortunately, the derivatives $f_\beta$ resp. $f_{\nu^\beta}$  grow quadratically resp. linear in  $\nnorm{\nu}$ contrary to the assumptions in \rp{6.1}.
Therefore, we define a smooth cut-off function $\tht\in C^\un(\R[]_+)$, $0<\tht\le2k$, where $k\ge k_0>1$ is to be determined later, by 
\begin{equation}\lae{6.10}
\tht(t)=
\begin{cases}
 t, &0\le t\le k,\\
2k , &2k\le t,
\end{cases}
\end{equation}
such that
\begin{equation}\lae{6.11}
0\le \dot \tht \le 4
\end{equation}
and consider the problem
\begin{equation}\lae{6.12}
\s_{2|_M}=f^{\frac12}(x,\tilde\nu),
\end{equation}
where for a spacelike hypersurface $M=\graph u$ with past directed normal vector $\nu$ we set
\begin{equation}\lae{6.13}
\tilde\nu=\tht(\tilde v)\tilde v^{-1}\nu.
\end{equation}
Then
\begin{equation}\lae{6.14}
\nnorm{\tilde \nu}\le ck,
\end{equation}
so that the assumptions in \rp{6.1} are certainly satisfied.
The constant $k_0$ should be so large that $\tilde \nu=\nu$ in case of the barriers $M_i, i=1,2.$
\rp{6.1} is therefore applicable leading to a solution $M_k=\graph u_k$ of \re{6.12}.

\nd From \ci[Lemma 8.1]{cg10} we then deduce that there exists a constant $m$ such that
\begin{equation}\lae{6.15}
\tilde v=(1-\abs{Du_k}^2)^{-\frac12}\le m \qquad \A k.
\end{equation}
Hence, $M_k=\graph u_k$ is a solution of \re{0.4}, if we choose $k\ge \max(2m, k_0)$.
\cvb

\vspace*{1.2cm}
\section*{Acknowledgement}

\cvb\nd
I would like to thank my advisor Prof. Dr. Claus Gerhardt for offering me the opportunity to work on this beautiful subject, for great support and for many interesting discussions.

\providecommand{\bysame}{\leavevmode\hbox to3em{\hrulefill}\thinspace}
\providecommand{\href}[2]{#2}

\end{document}